\documentclass[conference]{IEEEtran}
\IEEEoverridecommandlockouts
\usepackage{cite}
\usepackage{amsmath,amssymb,amsfonts}
\usepackage{algorithmic}
\usepackage{graphicx}
\usepackage{textcomp}
\usepackage{soul,color}
\usepackage{xcolor}
\usepackage{hyperref}

\newcommand{\R}{\mathbb{R}}
\newcommand{\T}{\mathsf{T}}
\renewcommand{\L}{\mathcal{L}}
\renewcommand{\b}{\boldsymbol}

\newcommand{\x}{\b{x}}

\newcommand{\lap}{\nabla^2}

\newcommand{\einf}{e_\infty}

\def\BibTeX{{\rm B\kern-.05em{\sc i\kern-.025em b}\kern-.08em
    T\kern-.1667em\lower.7ex\hbox{E}\kern-.125emX}}
\begin{document}

\title{A hybrid RBF-FD and WLS mesh-free strong-form approximation method\\
{\footnotesize
\thanks{
    The authors would like to acknowledge the financial support of the ARRS research core funding No.\ P2-0095, ARRS project funding No.\ J2-3048 and the World Federation of Scientists.}
}}

\author{
    \IEEEauthorblockN{Mitja Jančič}
    \IEEEauthorblockA{
        ``Jožef Stefan'' Institute \\
        Parallel and Distributed Systems Laboratory, \\
        Jamova cesta 39, 1000 Ljubljana, Slovenia \\
        and \\
        ``Jožef Stefan'' International Postgraduate School \\
        Jamova cesta 39, 1000 Ljubljana, Slovenia \\
        Email: mitja.jancic@ijs.si}
    \and
    \IEEEauthorblockN{Gregor Kosec}
    \IEEEauthorblockA{
        ``Jožef Stefan'' Institute \\
        Parallel and Distributed Systems Laboratory, \\
        Jamova cesta 39, 1000 Ljubljana, Slovenia \\
        Email: gregor.kosec@ijs.si}
}

\maketitle

\begin{abstract}
    Since the advent of mesh-free methods as a tool for the numerical analysis of systems of Partial Differential Equations (PDEs), many variants of differential operator approximation have been proposed. In this work, we propose a local mesh-free strong-form method that combines the stability of Radial Basis Function-Generated Finite Differences (RBF-FD) with the computational effectiveness of Diffuse Approximation Method (DAM), forming a so-called hybrid method. To demonstrate the advantages of a hybrid method, we evaluate its computational complexity and accuracy of the obtained numerical solution by solving a two-dimensional Poisson problem with an exponentially strong source in the computational domain. Finally, we employ the hybrid method to solve a three-dimensional Boussinesq's problem on an isotropic half-space and show that the implementation overhead can be justified.
\end{abstract}


\section{Introduction}
In recent years, mesh-free methods~\cite{wang2019methods} have been increasingly used to obtain a numerical solution to a system of PDEs. They are computationally more complex than traditional mesh-based methods, but the fact that they can operate on scattered nodes makes them very desirable, especially when complex three-dimensional domains are considered.

Since the advent of mesh-free methods in the 1970s, many different variants have been proposed, such as the Finite Point Method~\cite{onate1996finite}, the Generalized Finite Difference Method~\cite{gavete2003improvements}, the Diffuse Approximation Method (DAM)~\cite{prax1996diffuse} and the Radial Basis Function-Generated Finite Differences (RBF-FD)~\cite{tolstykh2003using}, to name but a few of the most commonly used, with recent research exploiting parallelism opportunities offered by a modern computer architecture~\cite{TrobecDepolli2021}.

While the RBF-FD is known for its high stability, DAM, also known as the Weighted Least Squares (WLS) approach, is known for its low computational complexity. Moreover, the WLS approach has been shown to be incredibly stable for low order approximations but has stability issues for higher order approximations~\cite{jancic2022stability}. On the contrary, the RBF-FD is stable even for higher order approximations. Thus, the aim of this paper is to combine the advantages of the RBF-FD variant (namely the stability) with the computationally efficient WLS variant by proposing a \emph{novel} hybrid WLS -- RBF-FD method. This method essentially splits the stencils into two separate sets: One that use the WLS approximation approach to approximate the differential operators and another one, that uses the RBF-FD approximation approach.

The stability and computational complexity of the proposed hybrid method are studied on a solution to a two-dimensional Poisson problem with an exponentially strong source~\cite{mitchell2013collection}. In addition, we also provide a solution to a three-dimensional Boussinesq's problem of the concentrated normal traction acting on an isotropic half-space~\cite{nwoji2017solution,slak2019adaptive}. We show that the hybrid method is more stable than the pure WLS variant and computationally cheaper than the pure RBF-FD variant.

\section{Solution procedure employing mesh-free methods}
\label{sec:solution_procedure}
\begin{figure}
    \centering
    \includegraphics[width=0.9\linewidth]{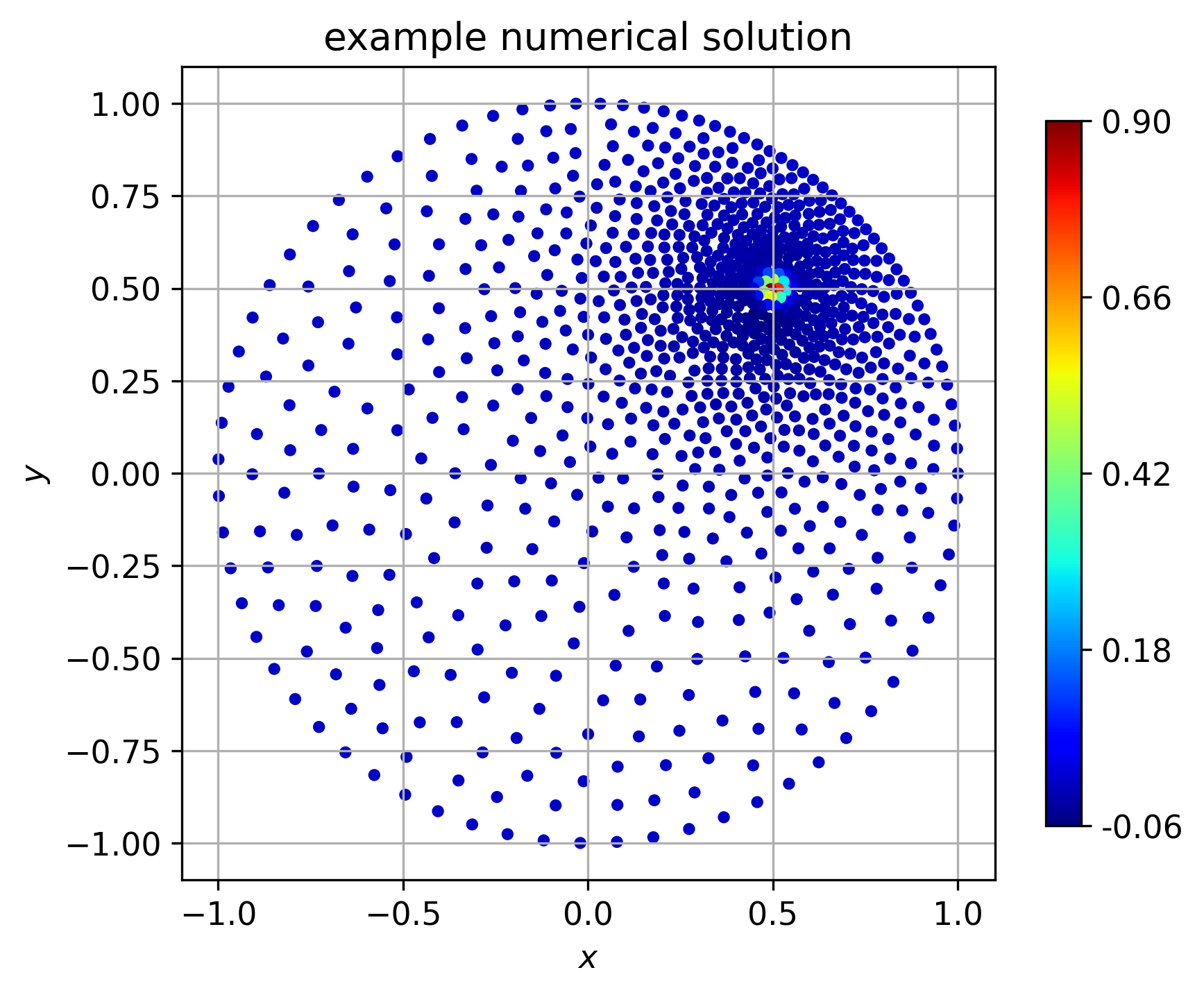}
    \caption{Example numerical solution on scattered nodes.}
    \label{fig:example}
\end{figure}

To obtain a numerical solution $\widehat{u}$ to a system of PDEs, three steps are required. First, the computational domain $\Omega$ is discretized using a dedicated node positioning algorithm that supports a spatially variable nodal distribution~\cite{slak2019generation} with a quasi-uniform internodal spacing $h$.
An example of nodal distribution is shown in Figure~\ref{fig:example}. A parallelized version of the same algorithm was recently published in~\cite{depolli_parallel_2022}, however, parallel execution is already out of the scope of this paper.

After discretizing the domain, the differential operators are approximated. A detailed procedure on differential operator approximation in the context of mesh-free methods is described in the following Section~\ref{sec:operator_approx}.

In the final step, the system of PDEs is discretized in spatial and temporal sense, resulting in a global system of linear equations. The system is solved and its solution is proclaimed as the numerical solution $\widehat{u}$ of a considered system of PDEs, of course, subject to given initial and boundary conditions.

\section{Linear differential operator approximation in the context of mesh-free methods}
\label{sec:operator_approx}
Consider a $d$-dimensional domain $\Omega \subset \R^d$ with a set of $N$ discretization points $\left\{ \x_j \right\}_{j=1}^N$. In the context of mesh-free methods, a linear differential operator $\L$ in any node $\x_c \in \Omega$ is approximated over a set of neighboring nodes
\begin{equation}
    \label{eq:ansatz}
    \widehat{\L u}(\b x_c) = \sum_{i = 1}^n w_iu(\b x_i)
\end{equation}
for an arbitrary function $u$, weights $w_i$ yet to be determined and support domain size $n$ also called \emph{stencil size}. It has been reported that a well-designed stencil can significantly reduce the computational cost~\cite{davydov2022improved}, but, usually, as is also the case in this paper, the closest $n$ nodes are chosen as the stencil to a central node $\x _c$.

The weights from equation~\eqref{eq:ansatz} are calculated for a given set of $s$ basis functions $\left\{ p _j \right\} _{j=1}^s$. In the Diffuse Approximation Method, a set of monomials $p_1, \dots,p_s$ with up to and including degree $m$ with $s=\binom{m+d}{d}$ are used as the approximation basis. This essentially means that the approximation~\eqref{eq:ansatz} can be compactly written as
\begin{equation}
    \label{eq:mon-system}
    \b{\mathrm{P}} \b w = \b\ell_p,
\end{equation}
where matrix $\b{\mathrm{P}}$ is a $n \times s$ matrix of monomials evaluated at stencil nodes,
\begin{equation}
    \b{\mathrm{P}} = \begin{bmatrix}
        p_1(\x_1) & \cdots & p_s(\x_1) \\
        \vdots    & \ddots & \vdots    \\
        p_1(\x_n) & \cdots & p_s(\x_n) \\
    \end{bmatrix}
\end{equation}
and $\b \ell_p$ is the vector of values assembled by applying the considered operator $\L$ to
the monomials at a central point $\x_c$
\begin{equation}
    \ell_p^i = (\L p_i(\x))\big|_{ \b x = \b x_c}.
\end{equation}

When the number of basis functions is equal to the stencil size ($s=n$), the described formulation yields a quadratic system of equations. Experience shows that the choice of larger support sizes ($n > s$) can be advantageous for accuracy and stability reasons. This leads to an overdetermined linear system, which is usually treated as a minimization of the Weighted Least Squares (WLS) norm. In the remainder of this paper, the above approximation approach, using only monomials as basis functions, will be referred to as the WLS approximation approach or WLS method.

Note that the same approximation procedure~\eqref{eq:ansatz} can also be used for any other commonly chosen basis functions, such as Multiquadrics, Gaussians, Radial Basis Functions (RBFs). In this paper we focus on two different types of basis: the already presented approach with monomials only and the approach using Polyharmonic Splines (PHS) augmented with monomials. The latter leads to a RBF-FD variant of the mesh-free methods described in the following section.

\subsection{The mesh-free RBF-FD variant}
\label{sec:rbffd}
We now take RBFs $\varphi(\x) = \varphi(\left\| \x - \x_c\right\|)$ centered at the stencil nodes of a central node $\b x_c$. The approximation~\eqref{eq:ansatz} then takes a compact form
\begin{equation}
    \label{eq:rbf-system}
    \b \Phi \b w  = \b\ell_\varphi
\end{equation}
for matrix $\b \Phi$ of evaluated radial basis functions
\begin{equation}
    \b \Phi = \begin{bmatrix}
        \varphi(\left\| \x_1 - \x_1\right\|) & \cdots & \varphi(\left\| \x_n - \x_1\right\|) \\
        \vdots                               & \ddots & \vdots                               \\
        \varphi(\left\| \x_1 - \x_n\right\|) & \cdots & \varphi(\left\| \x_n - \x_n\right\|) \\
    \end{bmatrix}
\end{equation}
and $\b \ell_\varphi$ is the vector of values assembled by applying the considered operator $\L$ to the RBFs at central point $\x_c$
\begin{equation}
    \ell_\varphi^i = (\L \varphi(\left\| \x - \x_i\right\|)\big|_{ \b x = \b x_c}.
\end{equation}

Different RBFs can be used. To avoid the dependency on a shape parameter, we choose Polyharmonic splines (PHS)
\begin{equation}
    \varphi(r) = \begin{cases}r^k,       & k \text{ odd}  \\
        r^k\log r, & k \text{ even}\end{cases},
\end{equation}
where $r$ denotes the Eucledian distance between two nodes. However, the approximation with a pure RBF basis guarantees neither convergent behavior nor solvability. To mitigate these problems, the approximation is augmented with a monomial basis by additionally enforcing an exactness constraint  for monomials, as we did in equation~\eqref{eq:mon-system}. This ensures convergent behavior and also allows us to control the order of the approximation, since the approximation order is the same as the order of the augmented monomials. This procedure finally results in a compactly written system
\begin{equation} \label{eq:rbf-system-aug}
    \begin{bmatrix}
        \b{\b \Phi}       & \b{\mathrm{P}} \\
        \b{\mathrm{P}}^\T & \b 0
    \end{bmatrix}
    \begin{bmatrix}
        \b w \\
        \b \lambda
    \end{bmatrix}
    =
    \begin{bmatrix}
        \b\ell_\varphi \\
        \b\ell_p
    \end{bmatrix}
\end{equation}
with Lagrangian multipliers $\b \lambda$. The system~\eqref{eq:rbf-system-aug} is overdetermined and treated as a constraint optimization problem~\cite{flyer2016role}. The weights are obtained by solving the system, while Lagrangian multipliers are discarded.

\subsection{Hybrid WLS--RBF-FD approximation approach}
\label{sec:hybrid}
The local RBF-FD systems~\eqref{eq:rbf-system-aug} are clearly larger than the purely monomial systems~\eqref{eq:mon-system}, making the RBF-FD method computationally more expensive. Therefore, our aim is to combine the computational efficiency of WLS approach with the high stability of RBF-FD variant to a create a novel \emph{hybrid} method.

The hybrid method has an additional step in the solution procedure, where we need to specify which stencils $\mathcal N (\x_i)$ use the WLS approach to approximate the differential operators and which the RBF-FD. This step essentially splits the $N$ discretization nodes of $\Omega$ into two parts: $N_{\text{WLS}}$ nodes whose stencils use WLS and $N_{\text{ RBF-FD }}$ nodes whose stencils use RBF-FD, where $N_{\text{ RBF-FD }} + N_{\text{WLS}} = N$.

Assigning a particular approximation type to a particular stencil is not a trivial task. The aim of a hybrid method is to ultimately result in numerical method that is more stable than the pure WLS and computationally less complex than the pure RBF-FD. Therefore, the RBF-FD approximation is naively assigned only to nodes with a high error of the numerical solution $\widehat{u}$ expectancy, while the rest are approximated with the WLS approach.

Normally, error indicators, such as ZZ-type~\cite{oanh2017adaptive}, are used in such cases. Although using an error indicator makes the most sense and would probably lead to better results, in this paper we make the decision a priori.

From an implementation point of view, only a small amount of overhead is required to implement a hybrid method. The biggest and practically only extra effort we have is when a global system is being assembled, because the sizes of WLS local systems~\eqref{eq:mon-system} and RBF-FD local systems~\eqref{eq:rbf-system-aug} do not match. Additional zero values have to be assigned in the global matrix to compensate for the mismatching sizes of the local approximations.

\subsection*{Note on the implementation}
All elements and corresponding functionality used in this paper are available as part of the \emph{Medusa library}~\cite{slak2021medusa}.

\section {Results}
In this section, an overview of the results is provided. We first study the proposed hybrid method on a two-dimensional Poisson problem with an exponentially strong source in the domain. In particular, we focus on the convergence rates and shape computation times. Finally, as a proof of concept, a three-dimensional Boussinesq's problem is solved in Section~\ref{sec:benchmark}.

All calculations were performed on a single core of a computer with
\text{Intel(R) Xeon(R) CPU E5-2620 v3 @ 2.40GHz} processor and 64 GB of DDR4 memory.
The code\footnote{Source code is available at \url{https://gitlab.com/e62Lab/public/cp-2022-splitech-hybrid-engine} under the tag \emph{v1.1}.} was compiled with \text{g++ (GCC) 9.3.0} for Linux with \text{-O3 -DNDEBUG} flags. The sparse system is solved using the single-threaded LU solver, unless otherwise specified.

\subsection{Two-dimensional synthetic example}
\label{sec:num_example}
The proposed hybrid method is studied by solving a synthetic example. We choose a two-dimensional elliptic PDE problem, i.e.\ a $d=2$ dimensional Poisson problem, with non-constant Dirichlet boundary conditions in domain $\Omega$. This example is usually used to test adaptive algorithms~\cite{mitchell2013collection}.

The problem is governed by
\begin{align} \lap u (\x) & = f_{\text{lap}}(\x) & \text{in } \Omega, \\ u (\x) & = f(\x) & \text{ on } \partial \Omega ,
\end{align}
where the domain $\Omega$ is a two-dimensional unit disc and the right-hand side is chosen to have an exponentially strong source
\begin{equation}
    f(\x) = \exp(-\alpha \left\| \x - \x_s \right\|^2),
\end{equation}
where $\alpha$ determines the strength of the source (for a strong source $\alpha =10^3$) and $\x_s = \b{1/2}$ is the location of the source. The Laplacian of $f(\x)$ can also be calculated analytically
\begin{equation}
    f_{\text{lap}} = 4(\alpha^2 \left\| \x - \x_s \right\|^2 - \alpha)\exp(-\alpha \left\| \x - \x_s \right\|^2).
\end{equation}

An example solution is shown in Figure~\ref{fig:example}. The above problem has an analytical solution $u(\x) = f(\x)$, which allows us to evaluate the accuracy of the numerically obtained solution $\widehat{u}$ in terms of the infinity norm error $\einf$.

The domain $\Omega$ was filled with $N$ scattered nodes with a variable node distribution that ensures the best local field description in the neighborhood of the strong source. In this work, the nodal distribution is given by
\begin{equation} h(\x)=\min(dx + (Dx-dx)\left\| \x - \x_s \right\|^{3/2}, dx),
\end{equation}
for $dx = Dx/5$ and 30 different values of $Dx$.

The problem was solved using all three previously described mesh-free variants, i.e.\ with the WLS approach using only monomials up to and including degree $m\in \left \{ 2, 4, 6 \right \}$, with the RBF-FD approach using Polyharmonic splines of order $k=5$ additionally augmented with monomials up to and including the same order $m$, and finally with a hybrid WLS--RBF-FD with the same approximation order. The stencil size $n$ was determined according to the recommendations of Bayona~\cite{bayona2017role} for a stable RBF-FD approximation
\begin{equation}
    n = 2\binom{m + d}{d}.
\end{equation}

\begin{figure}
    \centering
    \includegraphics[width=0.9\linewidth]{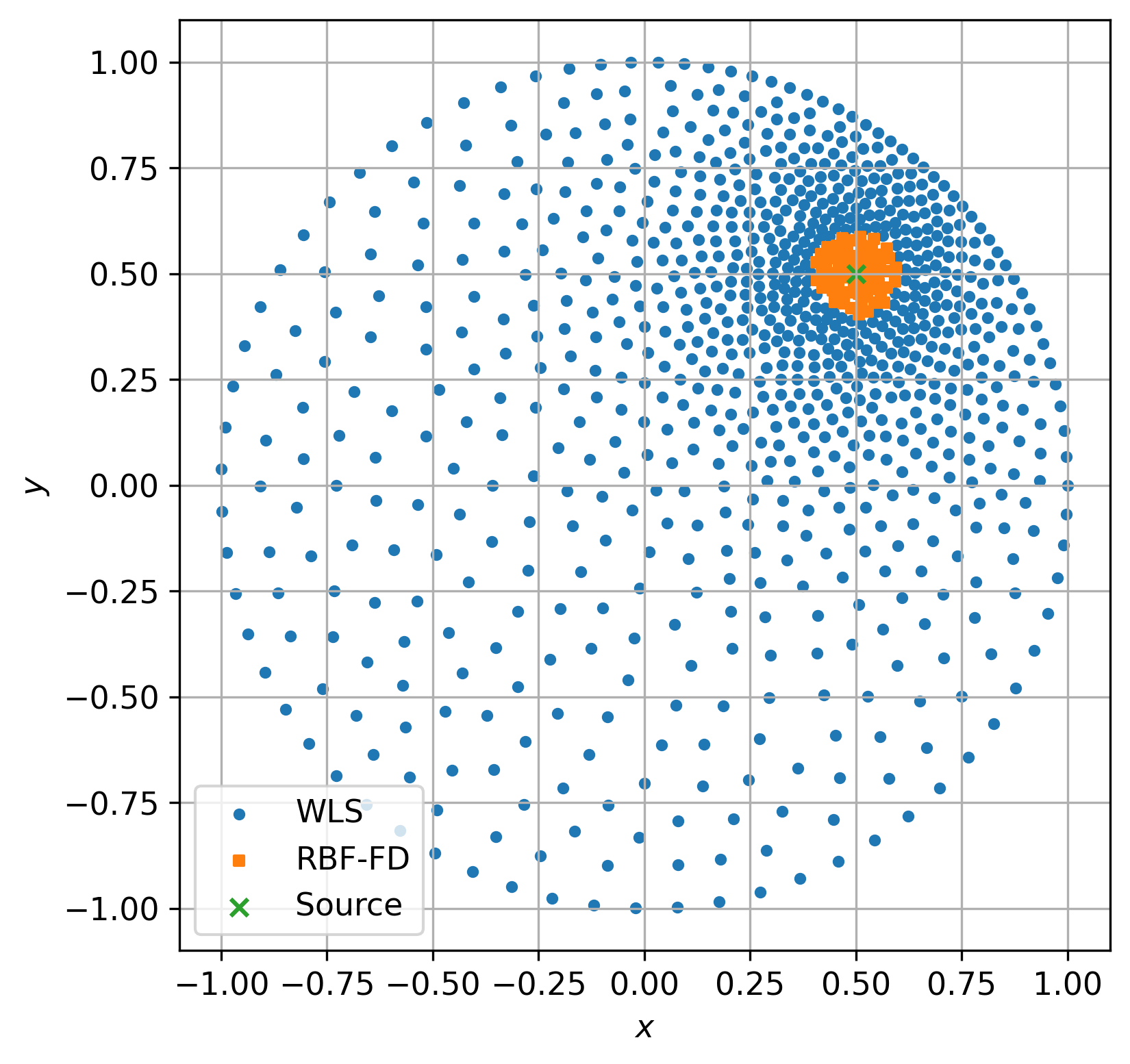}
    \caption{An example of different approximation methods used within the hybrid WLS--RBF-FD method.}
    \label{fig:approx-engines}
\end{figure}

The division of the nodes into $N_{\text{RBF-FD}}$ RBF-FD nodes and $N_{\text{WLS}}$ WLS nodes was done a priori without an error indicator. The largest error of the numerical solution is expected in the neighborhood of the exponentially strong source. We therefore define a circle with radius $r_s=0.15$ around the strong source $\x_s$. All the stencils with a central node $\x_c$ less than $r_s$ from the source are approximated using the more stable RBF-FD approach, while the rest use the WLS approximation. An example of the distribution of approximation types within the hybrid method is shown in Figure~\ref{fig:approx-engines}.

\subsubsection{Convergence rates}
In this paper, the error of the numerical solution is evaluated in computational nodes in terms of the infinity norm
\begin{equation}
    e_{\infty}  = \frac{\|\widehat{u} - u\|_\infty}{\|u\|_\infty}, \quad \|u\|_\infty = \max_{i=1, \ldots, N}
\end{equation}
because it shows the lowest convergence rates and, unlike the 2-norm, does not require averaging. After the solution is obtained, the computational nodes are discarded and the domain is discretized again with the same internodal spacing function $h$. With the new discretization, the shapes must be recomputed, which essentially allows us to evaluate the dependence of the approximation method on the quality of the node positioning. The process is repeated $N_{\text{runs}} = 100$ times, every time resulting in an approximately the same number of discretization nodes $N$.

\begin{figure*}
    \centering
    \includegraphics[width=\textwidth]{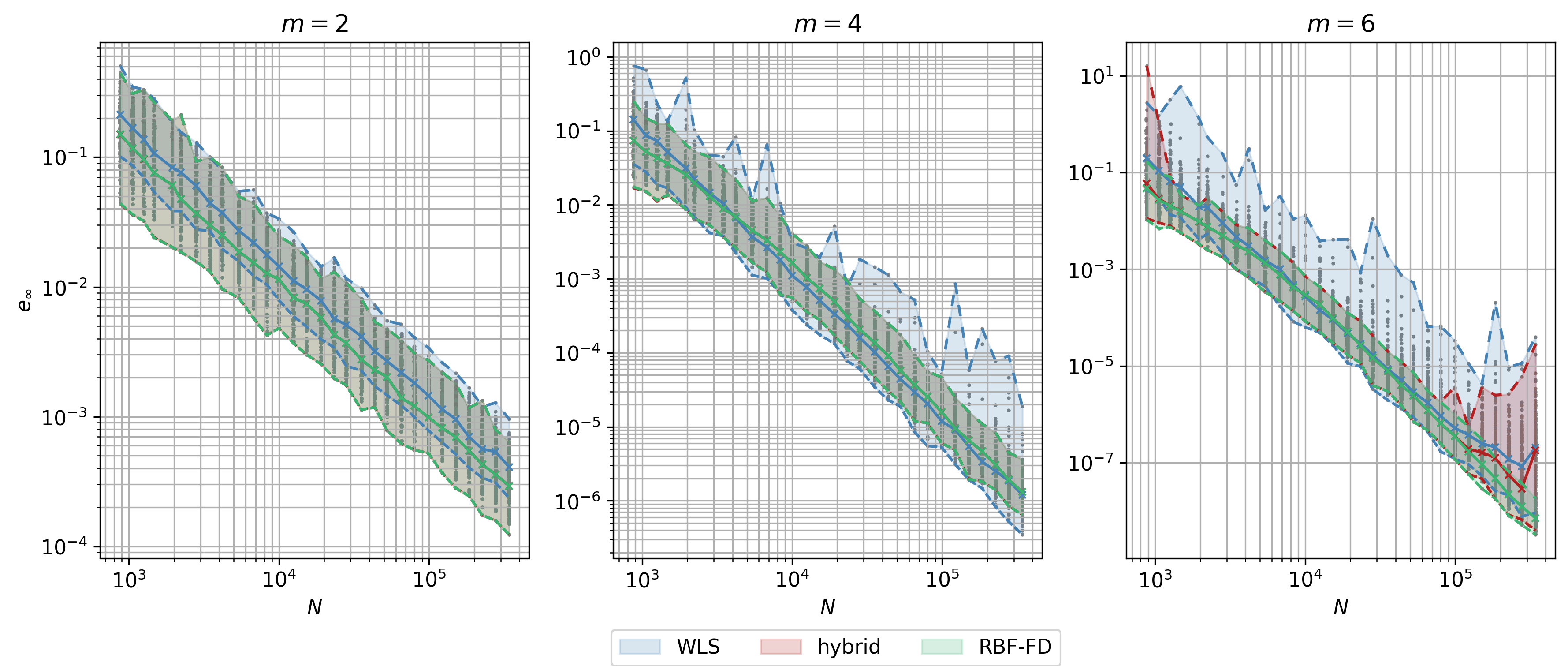}
    \caption{Convergence rates for WLS approximation approach (blue), RBF-FD approximation approach (green) and a novel hybrid approximation approach (red) for low order approximations $m=2$ (left) and higher order approximations $m=4$ (middle) and $m=6$ right.}
    \label{fig:convergences}
\end{figure*}

\begin{figure*}
    \centering
    \includegraphics[width=\textwidth]{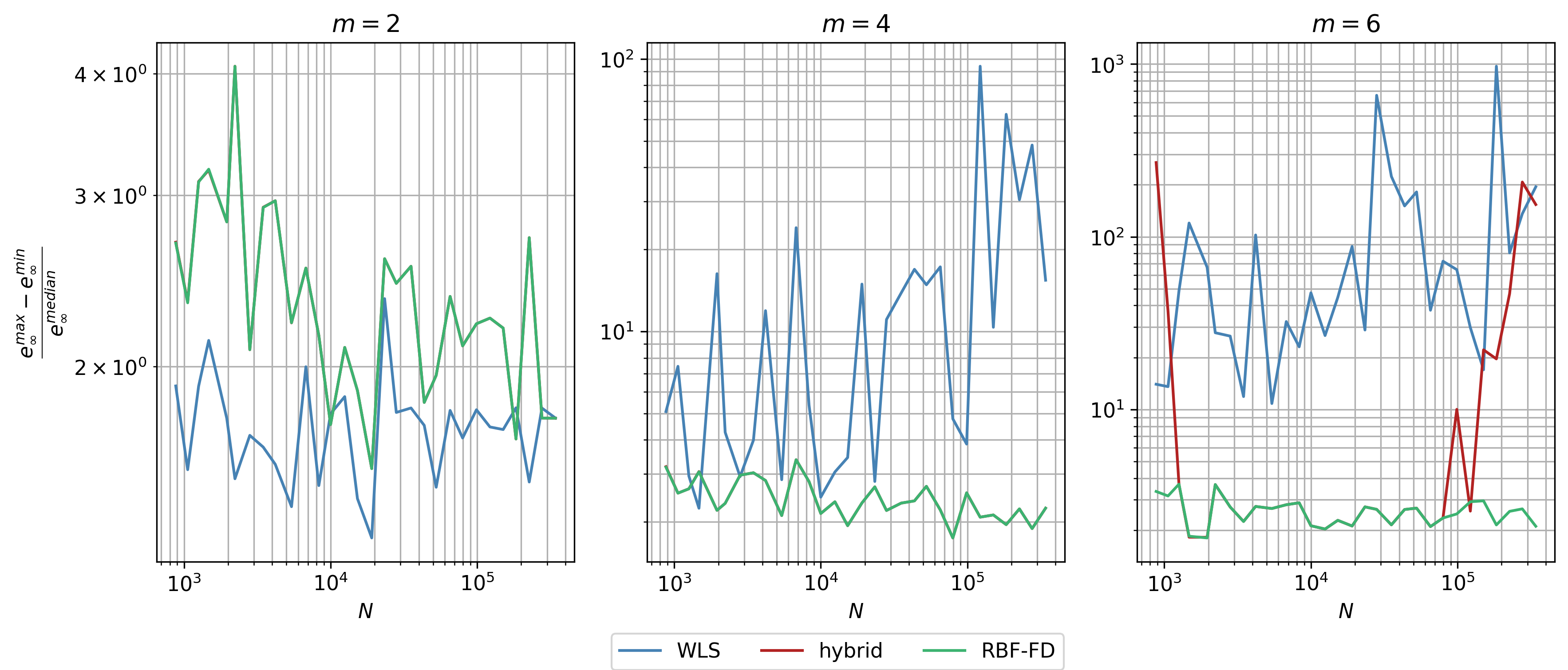}
    \caption{Normalized spread around a median infinity norm error after $N_{\text{runs}} = 100$ for WLS approximation approach (blue), RBF-FD approximation approach (green) and a novel hybrid approximation approach (red) for low order approximations $m=2$ (left) and higher order approximations $m=4$ (middle) and $m=6$ right.}
    \label{fig:spreads}
\end{figure*}

The convergence zones for three different approximation orders and three different mesh-free variants are shown in Figure~\ref{fig:convergences}. In the case of a low order approximation ($m=2$), we can immediately see that all three approximation methods are stable, with the smallest spread around a median of the infinity error norm clearly belonging to the WLS approach. The fact that the lower order WLS approximations are more stable with the WLS variant was already observed by Jan\v{c}i\v{c}~\cite{jancic2022stability}. Higher order approximations ($m=4$ and $m=6$), however, are more stable with the RBF-FD. The stability is further evaluated in Figure~\ref{fig:spreads}, making the advantages of a hybrid method in case of a higher order approximation even more evident.

\subsubsection{Computational times}
\begin{figure}
    \centering
    \includegraphics[width=\linewidth]{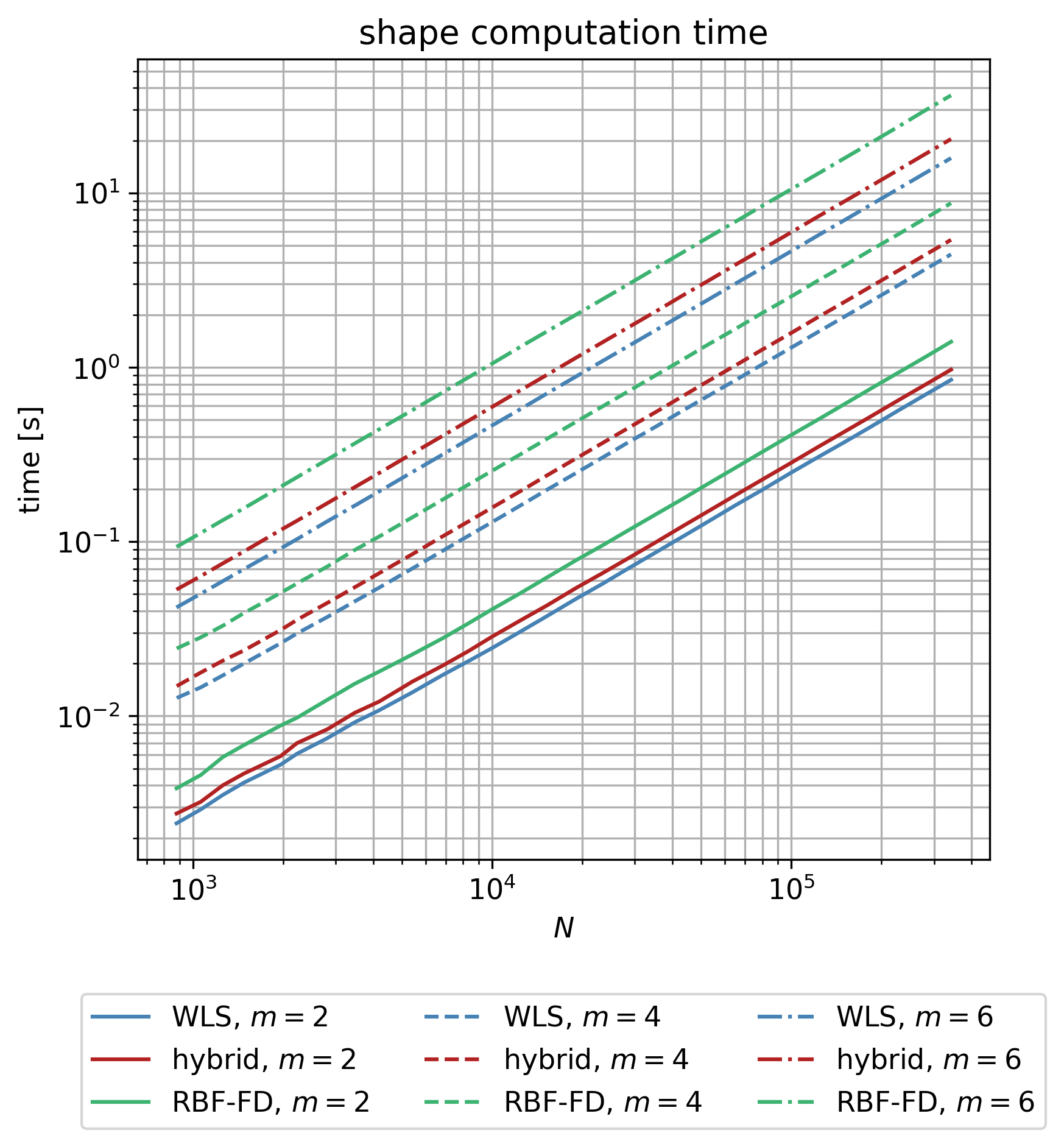}
    \caption{Shape computation times for WLS approximation approach (blue), RBF-FD approximation approach (green) and a hybrid WLS--RBF-FD approximation approach (red).}
    \label{fig:times}
\end{figure}

Another advantage of the hybrid method is that it is computationally cheaper than the pure RBF-FD approximation. This is demonstrated in Figure~\ref{fig:times}, showing the average shape calculation time for 10 runs. We can clearly see that the pure RBF-FD approximation is computationally expensive, while the WLS approach is computationally more efficient and the hybrid method is somewhere in-between - depending on the $N_{\text{RBF-FD}}/N$ ratio.

\subsection{Benchmark example}
\label{sec:benchmark}
As a benchmark case, we chose to solve a three-dimensional Boussinesq's problem of the concentrated normal traction acting on an isotropic half-space~\cite{nwoji2017solution,slak2019adaptive}. The problem is governed by the Cauchy-Navier equations
\begin{equation}
    (\lambda +\mu)\nabla (\nabla \cdot \b u)+\mu\nabla^2\b u = \b f
\end{equation}
with unknown displacement vector $\b u$, external body force $\b f$ and Lamé parameters
\begin{align}
    \lambda & = \frac{E \nu}{(1 - 2 \nu)/(1+\nu)} \text{ and } \\
    \mu     & = \frac{E}{2(1+\nu)},
\end{align}
for Young modulus $E=1$ and Poisson ratio $\nu=0.33$.

For domain $\Omega$ we take a three-dimensional box
\begin{equation}
    \Omega = \left \{ (x, y, z) \in \R^3,  -0.1 \leq (x, y, z) \geq -1 \right \}
\end{equation}
and discretize it using $h$-refinement towards the corner $\x_s = (-0.1, -0.1, -0.1)$ where force $\b P$ with magnitude 1 in the $-\widehat{e_z}$ direction is applied. The discretization resulted in a total of $N= 18849$ discretization points.

The problem has a closed form solution~\cite{nwoji2017solution} for displacements $\b u(\x) = \b u(x, y, z)$
\begin{align}
    u_x(\x) & = x\frac{P}{4\pi \mu}\Bigg(\frac{z}{\left\| \x \right\|^3} - \frac{1-2\nu}{\left\| \x \right\|(\left\| \x\right\|+z)}\Bigg) \\
    u_y(\x) & = y\frac{P}{4\pi \mu}\Bigg(\frac{z}{\left\| \x \right\|^3} - \frac{1-2\nu}{\left\| \x \right\|(\left\| \x\right\|+z)}\Bigg) \\
    u_z(\x) & = \frac{P}{4\pi \mu}\Bigg(\frac{z^2}{\left\| \x \right\|^3} + \frac{2(1-\nu)}{\left\| \x \right\|^3}\Bigg)
\end{align}
allowing us to calculate the infinity norm error in terms of the displacement magnitude.

To solve the sparse system, BiCGSTAB with ILUT preconditioner was used. The global tolerance was set to $10^{-14}$ with a maximum number of $500$ iterations, while the drop tolerance and fill factor were $10^{-5}$ and $30$ respectively.

\begin{figure}
    \centering
    \includegraphics[width=0.9\linewidth]{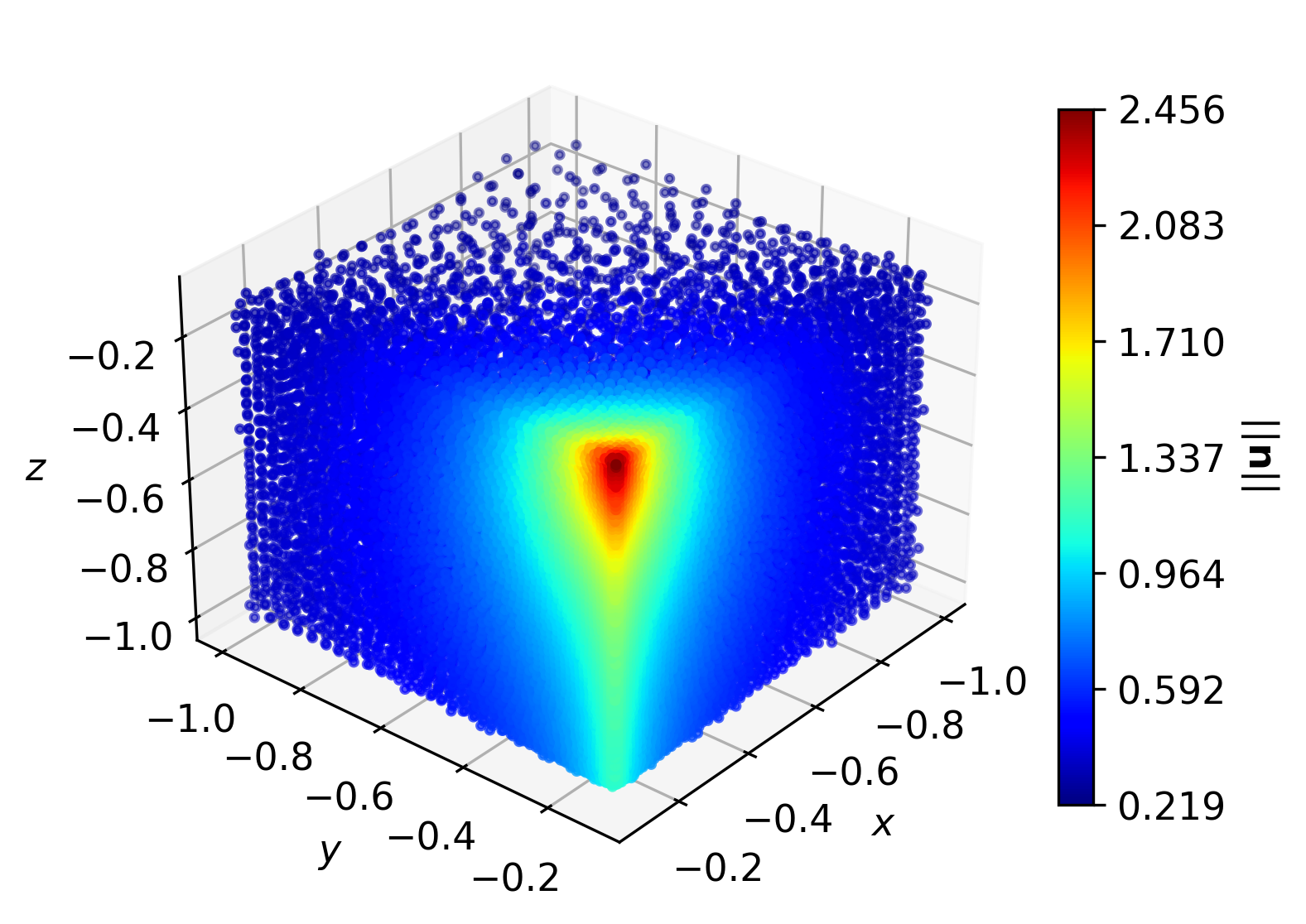}
    \caption{Example benchmark case. Numerical solution obtained with the proposed WLS--RBF-FD hybrid method.}
    \label{fig:contact}
\end{figure}

Results are computed using all three variants described previously, i.e.\ WLS (with Gaussian weights using $\sigma = 1.5$, essentially increasing the importance of nodes further away from the central stencil node), RBF-FD and a hybrid version of both with $r_s =0.5$, for monomials of order $m=4$ and PHS of order $k=5$. A visual representation of the solution obtained with the hybrid method is shown in Figure~\ref{fig:contact}, while a comparison of important numerical data is given in Table~\ref{tab:3d-data}.

\begin{table}
    \caption{Comparison table for the solution of  Boussinesq's problem.}
    \begin{center}
        \renewcommand\arraystretch{1.4}
        \begin{tabular}{cccc} \hline
            \multicolumn{1}{c}{Approximation} & \multicolumn{1}{c}{$\einf$} & \multicolumn{1}{c}{$t_{\text{shape}}$ [s]} & \multicolumn{1}{c}{$N_{\text{RBF-FD}} / N \cdot 100$}
            \\ \hline
            WLS                               & NaN                         & 4.74                                       & 0.00                                                  \\
            RBF-FD                            & 9.48$\cdot 10^{-5}$         & 8.22                                       & 100.00                                                \\
            hybrid                            & 2.37$\cdot 10^{-3}$         & 6.15                                       & 34.28                                                 \\ \hline
        \end{tabular}
        \label{tab:3d-data}
    \end{center}
\end{table}

We see that the novel hybrid method was able to obtain a numerical solution of sufficient quality. It is also clear from Table~\ref{tab:3d-data} that RBF-FD was able to achieve the best accuracy - approximately two orders of magnitude better than the hybrid method, but more importantly, the pure WLS approximation approach field to converge. This observation is of great importance, because it justifies the effort required to implement a hybrid method. It is also important to observe that less than 35 \% of the nodes from the hybrid method used the RBF-FD approximation approach, which is already enough to outperform the WLS in terms of stability and precision, and small enough to outperform the RBF-FD in terms of computation time, reducing it by about 33 \%.


\section{Conclusions}
\label{sec:conclusions}
A novel WLS--RBF-FD mesh-free method combining the RBF-FD and WLS variants is presented. We demonstrate that we can combine the advantages of the two commonly used mesh-free variants with only a small amount of additional work justified for the higher order ($m>2$) approximations.

Using a two-dimensional synthetic example with exponentially strong source, we show that the newly proposed hybrid method can be successfully used to obtain a numerical solution. We also demonstrate that the hybrid method is indeed computationally cheaper than the pure RBF-FD approach and more stable than the pure WLS approach for higher order approximations. Finally, on a solution to the three-dimensional Boussinesq's problem of the concentrated normal traction acting on an isotropic half-space we observe that the WLS variant fails to converge, while the hybrid WLS--RBF-FD method converges and reduces the shape computational times for about 33 \% compared to the pure RBF-FD.

In this work, the stencils were a priori divided into RBF-FD stencils and WLS stencils. We believe that better results could be obtained by using error indicators.

\bibliographystyle{IEEEtran}
\bibliography{references}

\end{document}